**Cognitive characteristics of intellectually gifted children with a diagnosis of ADHD.**

**Authors:** Cesare Cornoldi[1], David Giofrè[2], & Enrico Toffalini[1,**]

**Affiliations:**

[1]Department of General Psychology, University of Padua, Italy

[2]Disfor, University of Genoa, Italy

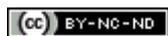

**Abstract**

Some children may be intellectually gifted, and yet experience behavioral and academic difficulties. We examined 82 twice exceptional children (2e-ADHD), having an excellent General Ability Index (GAI) derived from the Wechsler Intelligence Scale for Children-IV (GAI ≥ 125), and a diagnosis of Attention Deficit and Hyperactivity Disorder (ADHD). They accounted for 8.8% of a large sample of children with ADHD, which is twice as high as the proportion of intellectually gifted children in a typical population. This over-representation does not reflect a misdiagnosis of ADHD, as these children showed the typical features predicted on the grounds of data regarding the ADHD sample, including lower scores in working memory and processing speed measures, combined with the inclusion criteria for giftedness. Based on information concerning intellectually gifted children with either a Specific Learning Disorder (SLD) or typical development, we observed that these characteristics of intelligence are similar to those seen in SLD, but not in typical development, irrespective of whether 2e-ADHD children had a comorbid SLD.

**Cognitive characteristics of intellectually gifted children with a diagnosis of ADHD**

Attention problems associated with giftedness are very common, particularly among gifted children who often presents with some difficulties at school. A recent study reported that over 50% of children in a group of gifted under-achievers met the screening criteria for Attention Deficit and Hyperactivity Disorder (ADHD) based on teachers' reports, and almost one in three gifted under-achievers met the screening criteria for ADHD based on parents' reports (McCoach, Siegle, & Rubenstein, 2020). What is unclear, however, is whether gifted children's attention problems are due to an ADHD condition. Gifted children may easily become bored with the classroom activities typically proposed to children of their age, or be absorbed by their imagination and original ideas, without having ADHD. This situation raises doubts on the diagnosis of ADHD in gifted children (Hartnett, Nelson, & Rinn, 2004), and means that further research is warranted on the relationship between high intelligence and ADHD.

The aim of the present study was to collect new information on the cognitive characteristics of children who are considered twice exceptional ("2e"; Baum, Olenchak, & Owen, 1998) in that they are highly intelligent and have been concurrently diagnosed with ADHD. Such 2e-ADHD children have already been studied in terms of their behavioral and emotional problems (e.g., Foley-Nicpon, Rickels, Assouline, & Richards, 2012), but their cognitive characteristics have been under-investigated, and considered mainly by focusing on their attention deficits. For instance, Gomez, Stavropoulos, Vance, and Griffiths (2020) found that gifted children with ADHD tended to be more attentive than other children with ADHD, and Chae, Kim, and Noh (2003) reported that the former performed better than the latter in tasks measuring attention. It is therefore possible that gifted children receive a diagnosis of ADHD (Hartnett et al., 2004; Rinn & Reynolds, 2012) only because they present with poor attention in some circumstances, particularly when they become bored or are waiting for their classmates to complete a task. When Antshel et al. (2007) considered this problem, however, in a sample of 49 children diagnosed with ADHD who had an IQ higher than 120, they concluded that the diagnosis was appropriate because these children had other features typical of ADHD too (such as a family history of ADHD, and associated psychopathologies). Their study also included a few subtests of the WISC battery (Wechsler Intelligence Scale for Children) to examine the children's intellectual profile, but the resulting information was not used in discussing the identification and/or characterization of 2e-ADHD children.

In the present study we conducted a systematic analysis of the intellectual profile of a sufficient number of 2e-ADHD children (identified within a large sample of children with ADHD assessed for intelligence) in order to answer several questions that remain open. We did so by using the information available on a sample of around 1000 children with ADHD assessed with the most widely used intelligence test in the world, i.e., the WISC in its fourth version (Wechsler, 2004). This test battery provides information on various components of intelligence, grouped into four main areas, described by four corresponding Indexes, i.e., verbal (Verbal Comprehension Index, VCI) and non-verbal intelligence (Perceptual Reasoning Index, PRI) assessed with reasoning tasks based respectively on verbal and visuospatial materials, and two Indexes interesting respectively Working Memory (WMI), assessed with a simple and a complex span task requiring to remember alphanumerical material, and Processing Speed (PSI) assessed with tasks requiring to rapidly scan visual materials in order to give repeated appropriate responses. The test also provides an intelligence index, the General Ability Index (GAI), that is appropriate for use in cases of ADHD. In fact, the literature shows that some components of intelligence assessed by the WISC-IV, namely WMI and PSI, not only represent marked weaknesses in both ADHD and Specific Learning Disorders (SLD) (e.g., Kofler et al., 2018; Shanahan et al., 2006; Toffalini, Giofrè, & Cornoldi, 2017a; see also Ronald, de Bode, & Polderman, 2021, for an analysis of the genetic bases), but they are also more weakly related to general intelligence in both ADHD (Toffalini, Buono, & Cornoldi, 2022) and SLD (Giofrè & Cornoldi, 2015) than in the typical population. This outcome implies that it would be better to measure core intellectual abilities in cases of ADHD and SLD by means of the GAI. This index can be derived from the WISC-IV, by considering the measures of verbal and non-verbal reasoning, rather than the full-scale IQ (FSIQ; Saklofske, Prifitera, Weiss, Rolfhus, & Zhu, 2005). Other intelligence scales (e.g., Leiter-3; Roid, Miller, Pomplun, & Koch, 2013) could be used as they exclude WM, attention and PS separately from the general intelligence assessment too, but the WISC-IV has the advantage of considering a wider range of abilities, and it has a better consolidated psychometric reputation.

In the present study we examined four issues, as outlined below.

First, we sought to shed light on the extent and nature of the phenomenon, i.e., the proportion of intellectually gifted children within the ADHD population. If it is true that the association between intellectual giftedness and ADHD is relatively common, then the proportion of intellectually gifted children should be higher among children with ADHD than in the typically developing population, as already seen for children with SLD (Toffalini, Pezzuti, & Cornoldi, 2017). We also wanted to identify any specificities within the 2e-ADHD group, as regards gender distribution, age, and ADHD subtypes, with respect to the ADHD population as a whole. On the matter of gender, the typical greater proportion of males than females in ADHD populations (which seems particularly evident in Italy, see De Rossi et al., 2022) could even be emphasized in gifted children, as boys might be less motivated to use their intellectual abilities to compensate for their ADHD-related difficulties (Lai, Lin, & Ameis, 2022). Age seems to be relevant in ADHD (see Qian, Shuai, Chan, Qian, & Wang, 2013), and emerged as an important covariate in a study by Toffalini, Pezzuti, and Cornoldi (2017) on gifted children with SLD: with time these children became better able to cope with PSI tasks, but still had severe difficulties on WMI tasks. The picture becomes more complicated when we consider subtypes of ADHD as not all clinicians associate the diagnosis with a subtype, and there is some debate concerning whether different subtypes present different intellectual characteristics (e.g., Chhabildas, Pennington, & Willcutt, 2001; Fenollar-Cortés, Navarro-Soria, González-Gómez, & García-Sevilla, 2015; Mayes, Calhoun, Chase, Mink, & Stagg, 2009). Here too, examining specific populations of gifted children with ADHD could help to clarify this issue. In fact, Gomez et al. (2020) found that 2e-ADHD children had less difficulty with attentional aspects and greater difficulty with specific activity modulation and reflection comparatively to children with ADHD who were not intellectually gifted. They interpreted this result in the light of a systematic review conducted by Rommelse et al. (2015), who argued that intelligence may moderate the cognitive profile of individuals with neurodevelopmental disorders.

      A second aspect that our study examined concerns the average intellectual characteristics of 2e-ADHD children as a population, emerging when the four different aspects of intelligence assessed by the WISC-IV are separately considered, as suggested by the Manual itself (Wechsler, 2004). Relevant criticisms have been raised concerning the interpretation of the profiles of index scores of the WISC-IV at the individual level (e.g., Beaujean, 2017), also in the case of the Italian adaptation of the battery (Kush & Canivez, 2021). These criticisms are based on the consideration that most common variance across the subtests is explained by the "general" factor of intelligence, not by first-order factors. This does not exclude, however, that considerations at the population level might be important from both a clinical and scientific point of view. For example, Giofrè, Toffalini, Altoè, and Cornoldi (2017; see also Toffalini, Giofrè, & Cornoldi, 2017b) showed that mean differences between SLD and controls in some specific WISC-IV scores are so large that they potentially acquire a diagnostic significance, a result replicated with an even larger effect size by

Toffalini et al. (2022) on an ADHD sample. Also, Giofrè, Pastore, Cornoldi, and Toffalini (2019) showed that, unlike in the typical population, the first-order factors of intelligence effectively explain a larger portion of common variance than general intelligence in the SLD population. All these considerations could be particularly true in the case of children with ADHD who typically have good general reasoning abilities, but specific weaknesses in attention, PSI and WMI tasks (e.g., Kofler et al., 2018; Shanahan et al., 2006).

The third issue investigated here was whether information about a child's intellectual profile could help to allay any worries about a potential misdiagnosis of ADHD (Hartnett et al., 2004; Rinn & Reynolds, 2012). This matter could be considered in two ways: (i) by examining whether the intellectual profile of 2e-ADHD children resembles that of other children with ADHD, just on a higher overall level; and (ii) using a simulation procedure based on the consideration of the characteristics of the general population with ADHD. This latter type of procedure is becoming promising in the study of neurodevelopmental disorders (see Peters & Ansari, 2019). For example, recent studies (Carretti, Cornoldi, Antonello, Di Criscienzo, & Toffalini, 2022; Mammarella, Toffalini, Caviola, Lincoln, & Szucs, 2021) found that the mean profiles of domain-general and domain-specific cognitive abilities in children with either a reading or a math learning disorder could be inferred, with a good degree of approximation, from the set of linear correlations identified between the cognitive variables of interest on one side, and read and math ability scores on the other side, in the general population. Evidence is needed to show that the procedure can also be used in the case of giftedness. A study by Cornoldi, Giofrè, Mammarella, and Toffalini (2021), however, suggested that the characteristics of highly-gifted children did not match those predicted based on the general population. Specifically, the mean emotional response to testing in the population with very high achievement did not reflect the linear relationship observed between achievement score and emotional response in the general population. In the present study, a simulation procedure could be used to establish whether 2e-ADHD children match the characteristics that are expected to be found in the upper tail of the intelligence scores distributions considering the mean scores and covariances in the overall ADHD population.

The fourth matter that our study aimed to examine was whether 2e-ADHD children could also have SLD, as it is quite common to find an associated SLD in the ADHD population (see Pham & Riviere, 2015), and whether there could be differences between 2e-ADHD children with and without SLD. In fact, the few studies that distinguished between children who had ADHD with and without SLD generated different findings (probably due to differences in the samples' characteristics), arguing either for differences in their intellectual profiles (Becker, Daseking, & Koerner, 2021; Crisci, Caviola, Cardillo, & Mammarella, 2021; Katz, Brown, Roth, & Beers, 2011; Parke, Thaler, Etcoff, & Allen, 2020) or for surprising similarities (Toffalini et al., 2022). Examining cases of 2e-ADHD, identified on the strength of a high GAI, could shed important light on cognitive processes (WM and PS) related specifically to intellectual giftedness, by showing whether intellectual similarities between children with ADHD and those with SLD are due simply to the co-occurrence of the two disorders. It would be also worth examining whether 2e-ADHD children with SLD are cognitively weaker than 2e-ADHD children without SLD. As 2e-ADHD children with SLD are unable to compensate for their attentional failures and thus achieve sufficient learning levels, they would presumably have comparatively lower cognitive abilities than 2e-ADHD children without SLD. It could also be argued, however, that children with comorbid ADHD and SLD who perform well in intelligence tests have particularly high processing abilities to compensate for their difficulties.

To examine these four issues in the present study, we considered a sample of 2e-ADHD children, i.e., with a diagnosis of ADHD but also judged to be intellectually gifted on the grounds of a WISC-IV GAI, using a cut-off of 125, as done by other researchers (e.g.,

Assouline, Nicpon, & Whiteman, 2010). We first explored the size, gender and subtype distribution of this sample. Then we examined their intellectual performance in general, and specifically in measures not included in the calculation of the GAI, comparing them with two groups of children: one group of children with ADHD who were not intellectually gifted; and one group of children who were intellectually gifted and had not been diagnosed with ADHD. Where the necessary information was available, the gifted children with ADHD were divided into two subgroups, one with and the other without comorbid SLD. Finally, we simulated the intellectual profile of the 2e-ADHD group based on the assumption that these children represent one tail of the distribution of the group of children with ADHD as a whole, and we compared their predicted performance with their observed performance.

## 1. Method

### 1.1. Participants

The overall sample was the same as the one analyzed in a published study by Toffalini et al. (2022). It included 1,051 children (age range: 6.0 to 16.9 years, $M_{age}$ = 10.3, SD = 2.58; 80% males) diagnosed with ADHD collected by a network of selected clinicians working in different geographical areas (Northern vs Southern Italy) and in different types (private vs. public) of Centers of Italy. These clinicians are experts on ADHD, and share the same scientific background and assessment procedures, coordinated by a scientific association (AIRIPA) and the University of Padova, Italy. The clinicians provided anonymized data obtained when the WISC-IV was administered to the children under their supervision. The children had not been diagnosed with any other neurodevelopmental disorders, apart from SLD. They were native Italians and had been diagnosed with ADHD according to the criteria of the Diagnostic and Statistical Manual of Mental Disorders (5th ed.; DSM-5; American Psychiatric Association [APA], 2013) or the International

Classification of Diseases, Tenth Revision, Clinical (ICD-10-CM; World Health Organization, 1992). As typically done in the Italian Mental Health System, the ICD-10-CM codes for ADHD were adopted, sometimes only referring to a single, general category and code (F90), sometimes specifying one of the following four subtypes: F90.0, Attention deficit and hyperactivity disorder, predominantly inattentive type; F90.1, Attention deficit and hyperactivity disorder, predominantly hyperactive type; F90.2, Attention deficit and hyperactivity disorder, combined type; F90.8, Attention deficit and hyperactivity disorder, other type (typically used to specify conditions or terms like attention deficits in motor control and perception or undifferentiated attention deficit disorder). The information available regarding any pharmacological treatments, which are rarely administered in Italy (Bonati et al., 2021), was heterogeneous and unclear.

Clinicians were asked to specify whether they also had diagnosed any children at least 8 years old with SLD (in Italy this disorder is not routinely considered for younger children) based on the national guidelines for SLD. According to the Italian Consensus Conference on Learning Disorders (Istituto Superiore di Sanità, 2011), for a diagnosis of SLD a child's academic achievement should be consistently below the 5th percentile, or 2 SDs below average, in more than one test in at least one specific area of learning, in the absence of any socio-cultural or educational deprivation, or sensory, neurological, or intellectual impairments.

After omitting cases with any missing data regarding the WISC-IV indices or age, the sample consisted of 1004 children, which was reduced to 948 (age in [6.0, 16.8] years, Mage = 10.3, SD = 2.55; 81% males), when possible cases of intellectual disability (i.e., IQ <70, as in Toffalini et al., 2022) were removed. For the purposes of the present study, two subsamples of ADHD profiles were extracted and compared: one included 82 cases with a GAI of at least 125 (2e-ADHD); the other included 680 cases with an average GAI, of between 85 and 115 (average ADHD). See age and gender distributions below in Table 1

**Table 1**

*Comparison of the 2e-ADHD and average ADHD groups*

|  | 2e-ADHD | Average ADHD |
|---|---|---|
| Number of cases | 82 | 680 |
| as % of overall sample | 8.6% | 71.7% |
| Mean age (SD in brackets) | 10.13 (2.34) | 10.36 (2.59) |
| % of males | 87.8% | 80.4% |
| Mean FSIQ (SD in brackets) | 118.26 (8.21) | 93.17 (9.48) |
| Number of cases tested for SLD | 67 | 537 |
| Number of cases diagnosed with SLD (and % of the cases considered) | 35 (52.2%) | 275 (51.2%) |
| Number of cases with a specified subtype | 40 | 380 |
| *ICD-10 subtypes (N and % of the subtyped cases)* | | |
| Inattentive subtype (F90.0) | 19 (47.5%) | 120 (31.6%) |
| Hyperactive subtype (F90.1) | 0 (0.0%) | 7 (1.8%) |
| Combined subtype (F90.2) | 20 (50.0%) | 182 (47.9%) |
| Other (F90.8) | 1 (2.5%) | 71 (18.7%) |

*Note*. Data regarding SLD only concern children over 8 years old

For comparison, we also considered another sample of 129 intellectually gifted children with a GAI ≥125 and a diagnosis of SLD and no other neurodevelopmental disorders (age in [7.0, 15.8] years, Mage = 11.4, SD = 2.26; 69% males). This gifted SLD sample was selected from a larger database of 1622 WISC-IV profiles of children diagnosed with SLD who had no comorbidities (already analyzed by Toffalini, Giofrè, & Cornoldi, 2017a). The cases of SLD were diagnosed by specialist clinicians adopting the previously mentioned national guidelines, within the same AIRIPA network that provided the ADHD profiles.

**1.2. Instrument**

The Italian adaptation of the Fourth Edition of the Wechsler scale (WISC-IV) was used (Orsini, Pezzuti, & Picone, 2012). The Fifth Edition is not yet available in Italy, but the Fourth Edition seems particularly appropriate for studying developmental disorders (e.g., Mayes & Calhoun, 2007). The Italian version has internal consistencies, test-retest and inter-rater stability, and standard errors of measurement comparable with those of the English version (Wechsler, 2004). The children and adolescents were assessed individually in quiet rooms at the various centers involved in the study. The research was conducted in accordance with local institutional review board policies.

For the purposes of the present study, we examined the scores obtained in the ten core subtests of the WISC-IV, i.e., Block Design (BD), Similarities (SI), Digit Span (DS), Picture Concepts (PCn), Coding (CD), Vocabulary (VC), Letter-Number Sequencing (LN), Matrix Reasoning (MR), Comprehension (CO), and Symbol Search (SS). The children's raw scores were converted into the corresponding weighted scores (with $M = 10$, $SD = 3$) on the basis of the manual (Orsini et al., 2012), and of subsequently published instructions (Orsini & Pezzuti, 2014), to enable the scores of children of different ages to be compared. The weighted scores were used to calculate the IQ, the GAI, and the four factorial indexes described in the manual (Orsini et al., 2012; Wechsler, 2004). The IQ was calculated from the sum of the ten subtests, and the four main indexes from the sums of the relevant subtests, i.e.: the Perceptual Reasoning Index includes Block Design, Picture Concepts, and Matrix

Reasoning; the Verbal Comprehension Index includes Similarities, Vocabulary, and Comprehension; the Working Memory Index includes Digit Span and Letter-Number Sequencing; and the Processing Speed Index includes Coding and Symbol Search. We then calculated the scores for the GAI, obtained from the Verbal Comprehension and the Perceptual Reasoning Indexes (see also Flanagan & Kaufman, 2004; Wechsler, 2004).

**1.3. Analytical strategy**

The data analysis was conducted mainly on descriptive statistics, with inferences limited to binomial tests and mixed-effects linear models where relevant, as reported below (significance levels were obtained via likelihood ratio tests for nested models, based on the $\chi^2$). All data analysis was conducted with the R free software, version 4.2.1 (R Core Team, 2022).

A more complex aspect of the data analysis involved simulating mean profiles for gifted children with ADHD and typical development. These profiles were obtained via Monte Carlo simulation using a procedure similar to the one used in Giofrè et al. (2017) for typically developing children. For each group, we simulated a set of N observations reproducing the correlations between subtests identified in the respective populations (we used the large overall samples of ADHD [N = 948] and SLD [N = 1622] available to us, while for typical development we used the normative data from the WISC-IV Manual for the general population [N = 2200]) that were also obtained after removing cases of intellectual disability (Orsini et al., 2012). The means and SDs obtained in the subtests were also those observed in the respective populations (e.g., for typical development: M = 10, SD = 3 in all subtests). Then we selected simulated cases with a GAI ≥ 125 and recorded their simulated profiles. To account for uncertainty, we repeated this process for 10,000 iterations, and took median values as the central estimates for the means and SDs of the subtests; 95% CIs of the means were calculated from the standard errors derived from the simulated distributions.

## 2. Results

### 2.1. Numerosity and general characteristics of the 2e-ADHD sample

The 82 children who were 2e-ADHD accounted for 8.6% of the sample. Intelligence is assumed to follow a normal distribution, so only about 46 of the sample of 948 children with ADHD considered should have had a GAI ≥125 (i.e., 4.9% of the normal distribution after removing the lower tail of intellectual disability) if they simply reflected the characteristics of the overall distribution. Gifted children were over-represented in the ADHD sample even considering a GAI of 120 (15.0% vs 9.3% predictable based on the normal distribution) or 130 (3.4% vs 2.3%), and all these frequencies were significant based on a binomial test (all ps < 0.05). On the other hand, only 4.9%, 1.6% and 0.4% of the 948 children with ADHD had a full-scale IQ ≥120, or ≥ 125, or ≥ 130, respectively. Those percentages were considerably, and significantly (all ps < 0.001), below the percentage expected from the theoretical distribution based on a binomial test. Fig. 1 shows the whole distribution of the scores obtained on the GAI index in the larger ADHD sample (including those with IQ <70).

**Figure 1**

*Distribution of scores on the General Ability Index (GAI) in the ADHD sample. The black solid curve represents the expected normal distribution of the GAI scores in the general population (M = 100, SD = 15). Vertical dashed lines represent superimposed cutoffs.*

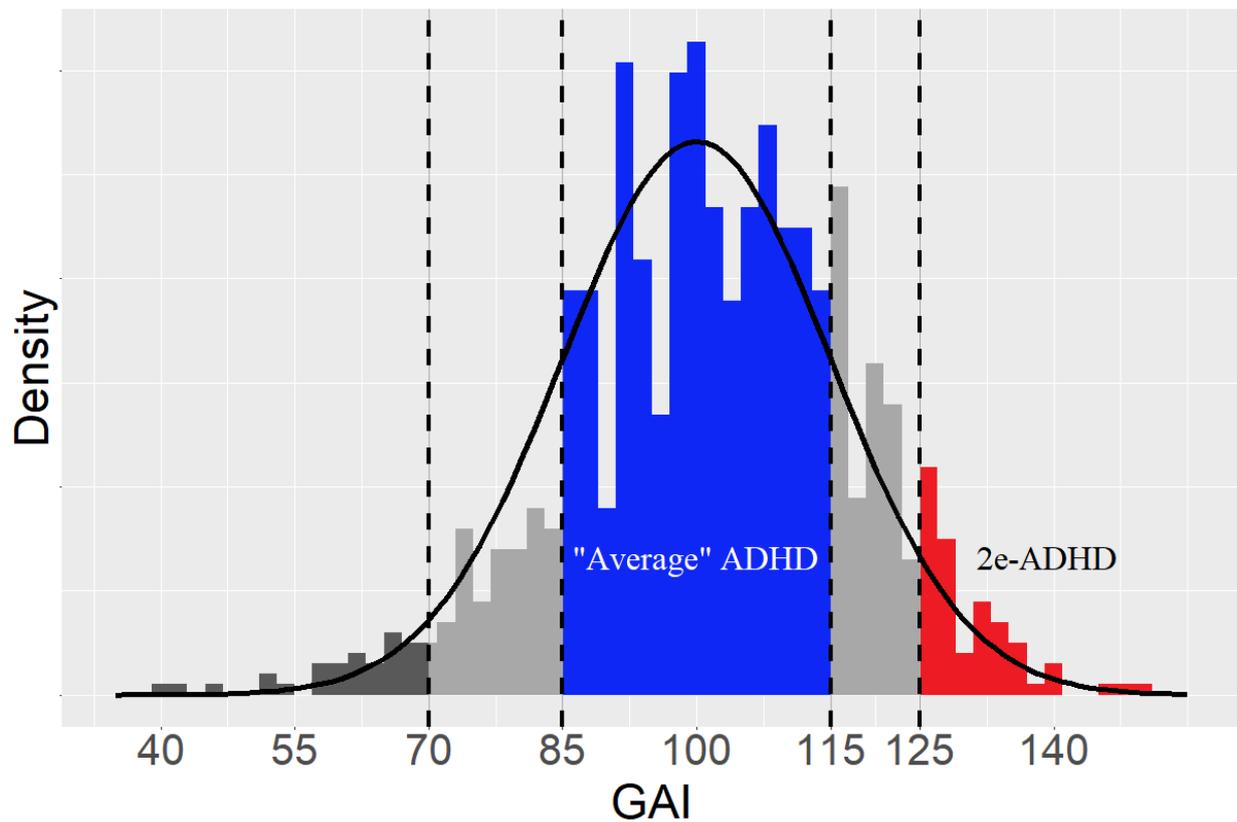

Focusing now on the group of 82 children who were 2e-ADHD, we first compared these children with a group of children with ADHD of average intelligence, i.e., with a GAI between 115 and 85, corresponding to one standard deviation above and below the mean (n = 680). Table 1 shows the characteristics of the two groups, which were of similar age. It is worth noting the similarity in the proportions of children with SLD (calculated only for the cases for which this information was available), which shows that the intellectual giftedness of children with 2e-ADHD does not lower their likelihood of having a SLD.

However, there were differences in the proportions of boys and ADHD subtypes in our samples. The proportion of boys was generally very high, as already reported for Italy (De Rossi et al., 2022), but boys also accounted for a particularly large proportion of the 2e-ADHD group (87.8% of cases). The proportion of the Inattentive subtype (F90.0) was larger in the 2e-ADHD group, whereas the Hyperactive (F90.1) and Other (F90.8) subtypes mainly concerned the average ADHD group.

### 2.2. Intellectual characteristics of the 2e-ADHD group

Table 2 shows the mean weighted scores obtained by the average ADHD group and the 2e-ADHD group in the ten basic WISC-IV subtests, and all the indices (first two columns of the table). In both groups, there was a marked discrepancy between the six subtests comprising the GAI and the four involving WM and PS. This discrepancy amounts to around 2.5 weighted points in the average ADHD group, and is even emphasized, approximately double (i.e., 1.7 SDs) in the 2e-ADHD group. The similarity between the profiles of the average ADHD group and the 2e-ADHD group supports the conviction that the latter were not misdiagnosed. Table 2 also shows the scores for the intellectually gifted (GAI ≥125) children with SLD (but not ADHD) extracted from the sample studied by Toffalini, Pezzuti, and Cornoldi (2017). The profiles of the 2e-ADHD and 2e-SLD children largely overlap.

Table 2. Mean WISC-IV scores (standard deviations in brackets) obtained in the ten basic subtests and in the derived Indexes by children with: a) ADHD and an average intelligence; b) ADHD and intellectual Giftedness; c) Specific Learning Disorder (SLD) and intellectual Giftedness.

**Table 2**

*Mean WISC-IV scores (standard deviations in brackets) obtained in the ten basic subtests and in the derived Indexes by children with: a) ADHD and an average intelligence; b) ADHD and intellectual Giftedness; c) Specific Learning Disorder (SLD) and intellectual Giftedness*

|      | Average ADHD (N = 675) | Gifted ADHD (N = 82) | Gifted SLD (N = 129) |
|------|------------------------|----------------------|----------------------|
| SI   | 9.79 (2.37)   | 14.49 (2.72)  | 14.91 (2.04)  |
| VC   | 9.99 (2.20)   | 14.41 (2.42)  | 14.32 (2.02)  |
| CO   | 10.19 (2.67)  | 14.68 (2.59)  | 14.94 (2.39)  |
| BD   | 10.03 (2.53)  | 13.60 (2.54)  | 13.48 (2.57)  |
| PCn  | 10.53 (2.71)  | 14.40 (2.31)  | 14.30 (2.15)  |
| MR   | 10.08 (2.48)  | 14.20 (2.28)  | 14.17 (2.25)  |
| DS   | 7.32 (2.54)   | 8.74 (2.56)   | 9.63 (2.80)   |
| LN   | 7.73 (2.40)   | 10.11 (2.27)  | 9.91 (2.26)   |
| CD   | 7.46 (2.93)   | 8.62 (2.92)   | 9.17 (3.13)   |
| SS   | 8.59 (2.87)   | 10.09 (3.11)  | 10.57 (2.86)  |
| VCI  | 99.97 (10.63) | 127.15 (10.80)| 128.47 (8.64) |
| PRI  | 101.19 (11.27)| 126.05 (9.72) | 125.82 (8.67) |
| WMI  | 85.19 (12.11) | 97.37 (11.54) | 98.60 (12.54) |
| PSI  | 88.19 (14.65) | 95.88 (15.24) | 98.84 (14.84) |
| GAI  | 100.63 (8.54) | 130.24 (5.58) | 130.60 (5.32) |
| FSIQ | 93.17 (9.48)  | 118.26 (8.21) | 120.43 (7.11) |

*Note*. Giftedness is defined as having a GAI of 125 or above. WISC-IV subtests: SI = Similarities, VC = Vocabulary; CO = Comprehension; BD = Block Design; PCn = Picture Concepts; MR = Matrix Reasoning; DS = Digit Span; LN = Letter-Number Sequencing; CD = Coding; SS = Symbol Search; WISC Indexes: VCI = Verbal Comprehension Index; PRI = Perceptual Reasoning Index; WMI = Working Memory Index; PSI = Processing Speed Index; GAI = General Ability Index; FSIQ = Full-Scale IQ.

### 2.3. Comparison between the 2e-ADHD children with and without SLD

A further comparison concerned only the 2e-ADHD group, distinguishing the children with ADHD alone from those with an associated diagnosis of SLD. The results are shown in Table 3. These findings should be taken with caution as they concern limited numbers of cases (N = 32 and 35, respectively), but no relevant differences emerged between the two subgroups in any of the subtests or indices. None of the comparisons reached statistical significance, and the effect sizes (Cohen's d) were small.

**Table 3**

*Mean WISC-IV scores (standard deviations in brackets) obtained by the 2e-ADHD subgroups without (N = 32) and with (N = 35) an associated Specific Learning Disorder (SLD).*

|      | 2e-ADHD (no comorbidity) | 2e-ADHD (with SLD) | t     | p   | Cohen's d |
|------|--------------------------|--------------------|-------|-----|-----------|
| SI   | 14.44 (2.47)             | 13.77 (2.67)       | 1.06  | .29 | 0.26      |
| VC   | 14.09 (2.53)             | 14.37 (2.22)       | -0.48 | .64 | -0.12     |
| CO   | 14.66 (2.81)             | 14.54 (2.69)       | 0.17  | .87 | 0.04      |
| BD   | 13.69 (2.36)             | 13.89 (2.93)       | -0.31 | .76 | -0.07     |
| PCn  | 14.72 (2.19)             | 14.20 (2.22)       | 0.96  | .34 | 0.24      |
| MR   | 14.25 (2.13)             | 14.31 (2.49)       | -0.11 | .91 | -0.03     |
| DS   | 9.03 (2.40)              | 8.49 (2.37)        | 0.93  | .35 | 0.23      |
| LN   | 9.94 (2.40)              | 10.03 (2.08)       | -0.17 | .87 | -0.04     |
| CD   | 8.59 (2.82)              | 8.17 (2.84)        | 0.61  | .54 | 0.15      |
| SS   | 10.75 (2.64)             | 10.06 (3.46)       | 0.93  | .36 | 0.22      |
| VCI  | 126.44 (11.41)           | 125.26 (10.19)     | 0.44  | .66 | 0.11      |
| PRI  | 126.75 (8.61)            | 126.63 (10.80)     | 0.05  | .96 | 0.01      |
| WMI  | 98.22 (11.45)            | 95.63 (10.60)      | 0.96  | .34 | 0.24      |
| PSI  | 97.25 (14.03)            | 94.80 (16.67)      | 0.65  | .52 | 0.16      |
| GAI  | 130.34 (6.51)            | 129.43 (3.98)      | 0.69  | .50 | 0.17      |
| FSIQ | 119.38 (6.49)            | 117.40 (7.12)      | 1.19  | .24 | 0.29      |

*Note.* See Table 1 note for the subtest and index acronyms.

**2.4. To what extent do the characteristics of 2e-ADHD children represent an extreme case of those seen in children with ADHD as a whole?**

In a final analysis we examined the degree to which the characteristics of the 2e-ADHD group (after applying the GAI ≥125 diagnostic cut-off on the continuum) represented an extreme case of those seen in children with ADHD generally. In other words, we wanted to see how closely our 2e-ADHD group reflected the general characteristics (correlations between subtests, mean scores, variability) of our larger ADHD population, as opposed to having characteristics unique to giftedness associated with ADHD. This analysis had two main aims: to sustain the view that the characteristics of "atypical" populations (including the intellectually gifted) substantially reflect those of the general population (after imposing suitable cut-offs); and to provide evidence to support the validity of the diagnosis of ADHD for the children in our 2e-ADHD group.

Fig. 2 shows the estimated means of the weighted scores in all the WISC-IV basic subtests, both observed (in the overall population with ADHD and in the 2e-ADHD group) and simulated (in the 2e-ADHD group, and for the gifted population in the normative sample). The extensive overlap between the observed and simulated samples of the 2e-ADHD group confirms that our gifted children with ADHD did represent the extreme portion of the distribution of children with ADHD (reflecting the same structure of covariances and mean cognitive deficits), and that they had been diagnosed appropriately. The profile of the 2e-ADHD group seems to reflect characteristics in between the average ADHD group and the (simulated) gifted typically developing group, showing that the dramatic drop in performance in the subtests involving WM and PS is largely attributable to the children's ADHD, and only partly a consequence of characteristics specific to gifted children. In fact, the gifted typically developing group was selected on the basis of the same GAI cut-off but showed only a minor decline on these indices.

**Figure 2**

*Estimated mean weighted scores in the ten WISC-IV subtests for the 2e-ADHD group (both empirically observed and simulated from the ADHD population; see text), for gifted typically developing (TD) children (simulated), and for the average ADHD group. Error bars represent 95% CIs.*

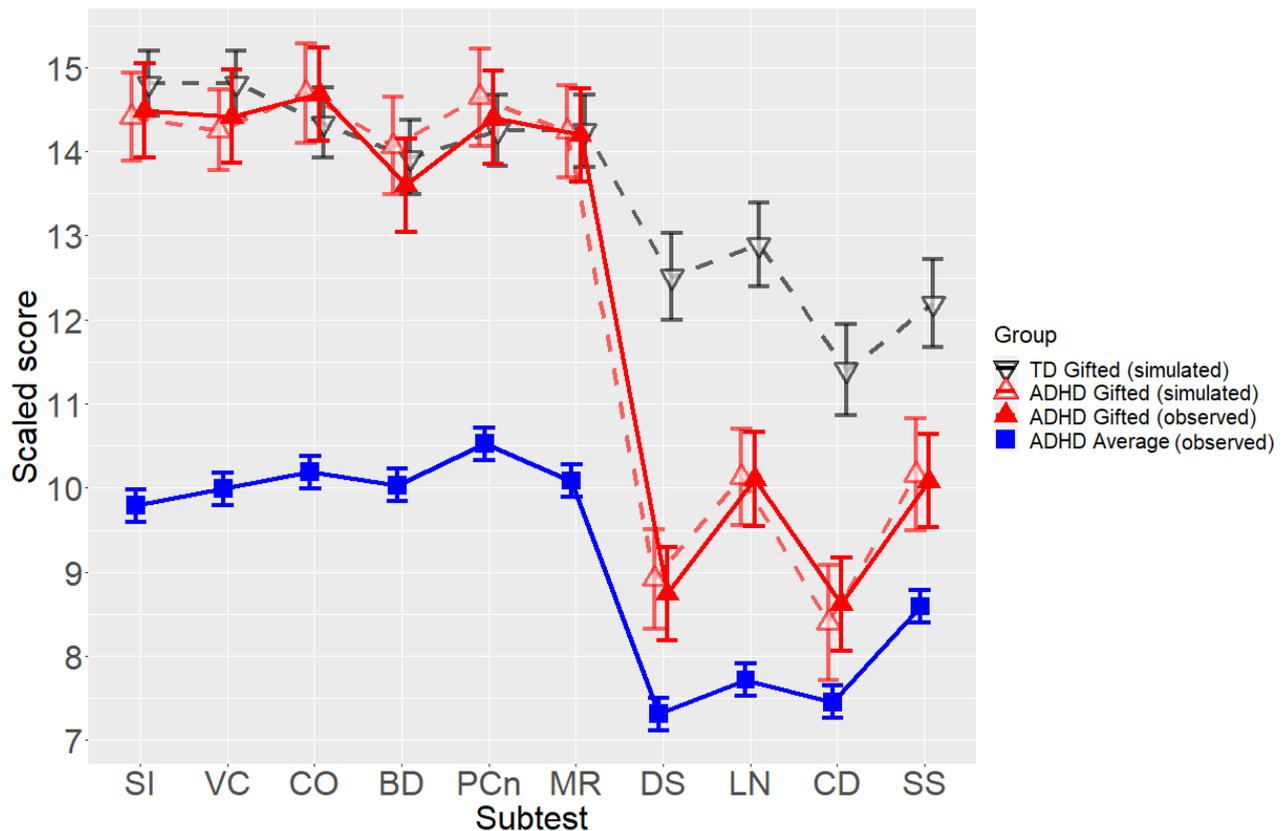

*Note*. See Table 1 note for the subtest and index acronyms.

## 3. Discussion

This study offers new and clearer information on the intellectual characteristics of 2e-ADHD children who are twice exceptional, as they are intellectually gifted and have ADHD. The study focused on four issues. The first concerned the general characteristics of our 2e-ADHD group. Although our sample was considerably large, our results should be considered with caution because a larger sample is probably needed to generalize regarding the characteristics the 2e-ADHD population. Our sample was drawn from a nationwide multicenter effort and is representative of the Italian population. It also reflects the diagnostic procedures shared by the centers involved, and typically adopted in Italy for the assessment of ADHD. For example, our sample included a larger preponderance of boys than in other studies (for a discussion on gender biases, see Garb, 2021) and this disproportion was even more evident in our 2e-ADHD group (possibly because gifted girls were better able to use their intellectual resources to control their ADHD symptoms). The study nonetheless generated useful information regarding the criteria for defining intellectual giftedness and the frequency of the association between a diagnosis of ADHD and intellectual giftedness. Using the traditional WISC criterion of a high full-scale IQ (including measures of WM and PS), cases of 2e-ADHD are much rarer than cases of giftedness in the typically developing population. The full-scale IQ seems inappropriate, however, because WM and PS measures reflect weaknesses specific to children who have ADHD, with loadings on the g-factor that are not particularly high (0.71 and 0.46, respectively; Toffalini et al., 2022). That is why we chose to exclude these measures from the estimation of general intelligence, as already done elsewhere in the literature (e.g., Roid et al., 2013; Toffalini, Pezzuti, & Cornoldi, 2017). This led us to find that intellectual giftedness is more common among children with ADHD than among children with typical development. This situation had already emerged in the case of SLD (Toffalini, Pezzuti, & Cornoldi, 2017), but is somewhat counterintuitive. Having greater intellectual resources might be expected to help children control their attentional difficulties, thereby reducing their likelihood of being diagnosed with ADHD, but this was not the case in the present study. It is also worth noting that the attentional subtype of ADHD was particularly in evidence in our 2e-ADHD group, which goes against what we might have expected based on previous reports (Gomez et al., 2020), and suggests that intellectually gifted children with ADHD are better able to control their hyperactivity than their inattention.

The genetic research may offer important elements for the consideration of the relationship between ADHD, cognitive abilities and intelligence. For example, Ronald et al. (2021) suggest that ADHD is associated with an overall polygenic risk score that also involves neuropsychological functioning. Other evidence suggests that there is a genetic basis for a negative association between IQ and ADHD (Kuntsi et al., 2004). The latter study, however, was performed on a large sample of typically developing children of five years of age at risk for a series of problems, and tested using only two tests of the WPPSI, which might limit IQ reliability. Jepsen, Fagerlund, and Mortensen (2009), who reviewed studies on the association between IQ and ADHD, also suggested a negative relationship, albeit rather modest. In fact, in our initial unselected sample (N = 1051), the average full-scale IQ of children with ADHD was clearly lower as compared to the general population (i.e., 93.7), and only the GAI presented an average slightly above 100 (i.e., 101.7). This finding is not particularly surprising, as it is similar to results obtained with Italian children with SLD (see Giofrè & Cornoldi, 2015). Nonetheless, Kuntsi et al. (2004) used a two-subtest approximation of the IQ that could be considered a very short form of the GAI. To speculate, there might be high heterogeneity in the ADHD population, and the genetic influence could be non-linear in the population, being weaker at high levels of intellectual functioning than it is in the general population.

The overrepresentation of intellectually gifted children in the ADHD overall sample considered in the present study might also be due to an underrepresentation of children with low intelligence scores. That is, families might be more likely to seek for a clinical diagnosis when poor behavioral regulation was judged as unexpected (and more malleable) considering the overall cognitive functioning. In addition, more affluent families could be more likely to seek clinical assessment. However, it is also possible that parents may look for the clinician's support when the situation has become serious, the child has a disorder that cannot compensate with general intellectual resources, and family resources are insufficient. In fact, in our study, different clinical centers were involved, some of which are currently supported by the Italian public health service and do not require extra fees for the family. In brief, our results should be corroborated by larger population-based epidemiological studies.

It is also worth mentioning that the Italian school context is different from other countries. For example, children in Italy are enrolled in normal schools. Once the diagnosis is performed, children with special needs (e.g., ADHD, SLD) may receive extra care and attention in their classes (for example, they are allowed to have extra time during tests and other compensatory measures). Therefore, we cannot rule out the possibility that our sample presents with some differences as compared to other international samples.

We then separately considered different intellectual aspects of children with ADHD, taking a 'strengths and weaknesses' approach. This approach has met with some criticism when used to draw conclusions at the individual level (Beaujean, 2017), but it appears informative when used to compare between populations, as shown in previous research (Giofrè et al., 2017; Toffalini, Giofrè and Cornoldi, 2017a, Toffalini, Giofrè and Cornoldi, 2017b), and also resulted useful in the present study. It confirmed that even intellectually gifted children with ADHD have difficulties relating to WM (e.g., Kofler et al., 2018) and PS (e.g., Shanahan et al., 2006). Though they perform better than children with ADHD of average intelligence, they retain an impressive imbalance between their scores in the GAI subtests and those concerning WM and PS. This imbalance cannot be due simply to the high GAI cut-off adopted, as applying the same selection criterion to the typically-developing population revealed a much smaller imbalance. It seems instead to reflect the combined effect of two phenomena seen in the ADHD population as a whole, namely the discrepancies in the intellectual profile, and the lower loadings of WM and PS on the g factor by comparison with the typically-developing population (cf. Toffalini et al., 2022).

When Toffalini, Pezzuti, and Cornoldi (2017) used a similar approach to find intellectually gifted children in a SLD population, a similar result emerged, but the pattern was more clearly influenced by age, showing that PS difficulties gradually tended to disappear, and WM difficulties tended to increase. Conversely, our analysis on the present sample of 2e-ADHD showed no clear age-related effects, although it suggested that WM difficulties (in the LN subtest) may increase with age, and difficulties with the WISC-IV PS tasks persist across ages. This last result can be interpreted bearing in mind that PS tasks place a strong demand in terms of maintaining attention, and children with ADHD might be quick to give a single response, but become slow and discontinuous when they need to stay focused in order to provide repeated responses (Borella, De Ribaupierre, Cornoldi, & Chicherio, 2013).

Examining the intellectual profile of our 2e-ADHD children generated relevant information concerning whether gifted children might be misdiagnosed with ADHD (Hartnett et al., 2004; Rinn & Reynolds, 2012) because they get bored, divert their attention elsewhere if a lesson is not challenging enough, or respond inappropriately to adults' reactions to their extreme precociousness. The similarities in the intellectual patterns of our average ADHD and 2e-ADHD groups suggests that the latter had been correctly diagnosed. Furthermore, the impressive overlap found between the data predicted based on the performances of the overall population of children with ADHD and the actual observations for the 2e-ADHD children also supports this conclusion. On this respect it should be noticed that the method of inferring the characteristics of a tail of the distribution from the overall population (see Carretti et al., 2022; Mammarella et al., 2021; Peters & Ansari, 2019) resulted fruitful, showing that the method can also be used in the case of the higher tail, i.e., of giftedness, despite the difficulties emerged when emotional responses to testing were considered (Cornoldi et al., 2021). Our findings, however, do not rule out the risk of a few gifted children being misdiagnosed, and clinicians should be cautious about making a differential analysis between ADHD and giftedness. They can consider factors such as whether a child's behavior is problematic in all or only some settings, whether their attention span or performance improves when they engage in more interesting activities, and how much their performance in a task varies (Webb & Latimer, 1993).

In general, it seems that our 2e-ADHD children had been correctly diagnosed, and the question is why they are so frequent. An often-mentioned hypothesis is that gifted children might typically present psychomotor overexcitability (Karpinski, Kolb, Tetreault, & Borowski, 2018; Piechowski & Colangelo, 1984), and high levels of stress and emotionality (Gaesser, 2018). The overexcitability hypothesis has been criticized (Vuyk, Krieshok, & Kerr, 2016), however, and there is evidence (Francis, Hawes, & Abbott, 2016) to suggest that the mental health of intellectually gifted individuals is better than that of the typical population. Another possibility is that children who have remarkable reasoning abilities but whose other psychological characteristics are average have an imbalanced intellectual profile that they struggle to manage (see Baum et al., 1998) or that leads them to dysfunctional behaviors. In particular the 2e-ADHD could be influenced by the very strong imbalance between much above-average mainly "top-down" abilities (e.g., fluid reasoning and verbal skills) and much lower mainly "bottom-up" basic cognitive processes (e.g., WM and PS). These children, for example, could adopt strategies and attitudes, based on their abilities and on top-down processes, that from their point of view are more rapid and functional, but do not meet the specific requests posed by the context.

Our study also investigated whether any of our 2e-ADHD children had comorbid SLD as this is common in the ADHD population (see Pham & Riviere, 2015), and what intellectual characteristics distinguished them from the rest of the 2e-ADHD group. The issue of SLD is controversial in the case of all children with ADHD. Some studies have argued for important differences between children who have ADHD with versus without SLD (Becker et al., 2021; Crisci et al., 2021; Katz et al., 2011; Parke et al., 2020), but a recent study found that this was not the case (Toffalini et al., 2022). Examining the case of twice-exceptional children (identified on the basis of a high GAI) brought more light on the cognitive processes that support academic learning (working memory and processing speed) specifically in the case of intellectual giftedness. It has been claimed that the co-occurrence of impairments in different areas of learning would have an additive effect on the burden on the individual's intellectual potential (Landerl, Fussenegger, Moll, & Willburger, 2009). This might mean that 2e-ADHD children who are unable to compensate for their attentional failures, and are consequently diagnosed with SLD, would have comparatively weaker cognitive abilities than 2e-ADHD children without SLD. This was not the case in our sample, possibly because children with learning difficulties who perform well in intelligence tests are generally very smart.

In our study, the simulated profile obtained using the correlations matrix and current results of the 2e-ADHD group were closely aligned. It is worth noting, however, that other studies found that at the very top of the giftedness performance relationships seem to become non-linear (Humphreys, Lubinski, & Yao, 1993). In the present study, we did not find evidence of any non-linear relationships in the current sample, as the performance of the simulated and real groups were very closed. It should be noted, however, that the Wechsler scales are explicitly designed to measure performance around the normal range of the intelligence distribution, and they are probably not very suitable for exceptional performances, largely exceeding three standard deviations above average, for example. For these reasons, we cannot rule out the possibility that some nonlinear relationships might exist at very top of the distribution. In fact, our results should be replicated using a much larger sample size and using a battery explicitly developed for testing extremely gifted children.

In conclusion, our study offers a series of important indications that have theoretical, empirical and clinical implications. In particular, it confirms that it is crucial to examine the intellectual profile of children diagnosed with ADHD, and this may help to identify and support children who not only have ADHD, but are also intellectually gifted. This has important educational implications as teachers should adopt a personalized approach to such children, taking advantage of their strengths (e.g., inviting them to help other children, take on particular responsibilities, focus on their interests), and also considering their weaknesses (e.g., suggesting strategies for managing working memory and processing speed deficits, recognizing academic difficulties, avoiding excessive time pressures).

INTELLECTUAL GIFTEDNESS AND ADHD                                                                30Fenollar-Cortés, J., Navarro-Soria, I., González-Gómez, C., & a García-Sevilla, J. (2015). Detección de perfiles cognitivos mediante WISC-IV en niños diagnosticados de TDAH: Existen diferencias entre subtipos?. *Revista de Psicodidáctica, 20*(1), 157-176. https://doi.org/10.1387/RevPsicodidact.12531

Flanagan, D. P., & Kaufman, S. (2004). *Essentials of Assessment with WISC-IV*. New York, NY: Wiley.

Foley-Nicpon, M., Rickels, H., Assouline, S. G., & Richards, A. (2012). Self-esteem and self-concept examination among gifted students with ADHD. *Journal for the Education of the Gifted, 35*(3), 220-240. https://doi.org/10.1177/0162353212451735

Francis, R., Hawes, D. J, & Abbott, M. (2016). Intellectual giftedness and psychopathology in children and adolescents: A systematic literature review. *Exceptional Children, 82*(3), 279–302. https://doi.org/10.1177/0014402915598779

Gaesser, A. H. (2018). Befriending anxiety to reach potential: Strategies to empower our gifted youth. *Gifted Child Today, 41*(4), 186–195. https://doi.org/10.1177/1076217518786983

Garb, H., N. (2021). Race bias and gender bias in the diagnosis of psychological disorders. *Clinical Psychology Review, 90*:102087. https://doi.org/10.1016/j.cpr.2021.102087

Giofrè, D., & Cornoldi, C. (2015). The structure of intelligence in children with specific learning disabilities is different as compared to typically development children. *Intelligence, 52*, 36–43. https://doi.org/10.1016/j.intell.2015.07.002

Giofrè, D., Toffalini, E., Altoè, G., & Cornoldi, C. (2017). Intelligence measures as diagnostic tools for children with specific learning disabilities. *Intelligence, 61*, 140–145. https://doi.org/10.1016/j.intell.2017.01.014

INTELLECTUAL GIFTEDNESS AND ADHD                                              33oppositional-defiant disorder. *Child Neuropsychology, 13*:6, 469-493. https://doi.org/10.1080/09297040601112773

Mayes, S. D., Calhoun, S.L., Chase, G.A., Mink, D.M., & Stagg, R.E. (2009). ADHD subtypes and co-occurring anxiety, depression, and oppositional-defiant disorder: differences in Gordon diagnostic system and Wechsler working memory and processing speed index scores. *Journal of Attention Disorders, 12*(6), 540-550. https://doi.org/10.1177/1087054708320402

McCall, R. B. (1977). Childhood IQ's as predictors of adult educational and occupational status. Science, 197(4302), 482–483. https://doi.org/10.1126/science.197.4302.482

McCoach, D. B., Siegle, D., & Rubenstein, L. D. (2020). Pay attention to inattention: Exploring ADHD symptoms in a sample of underachieving gifted students. *Gifted Child Quarterly, 64*(2), 100-116. https://doi.org/10.1177/0016986219901320

Orsini, A., & Pezzuti, L. (2014). L' indice di abilità generale della scala WISC-IV [The WISC-IV General Ability Index]. *Psicologia Clinica Dello Sviluppo, 18*(2), 301–310. https://doi.org/10.1449/77640

Orsini, A., Pezzuti, L., & Picone, L. (2012). WISC-IV: *Contributo alla taratura Italiana. [WISC-IV Italian Edition]*. Florence, Italy: Giunti O. S.

Parke, E. M., Thaler, N. S., Etcoff, L. M., & Allen, D. N. (2020). Intellectual profiles in children with ADHD and comorbid learning and motor disorders. *Journal of Attention Disorders, 24*(9), 1227-1236. https://doi.org/10.1177/1087054715576343

Peters, L., & Ansari, D. (2019). Are specific learning disorders truly specific, and are they disorders? *Trends in Neuroscience and Education, 17*:100115. https://doi.org/10.1016/j.tine.2019.100115

INTELLECTUAL GIFTEDNESS AND ADHD                                                    34Pham, A. V., & Riviere, A. (2015). Specific learning disorders and ADHD: Current issues in diagnosis across clinical and educational settings. *Current Psychiatry Reports, 17*:38. https://doi.org/10.1007/s11920-015-0584-ygaesser

Piechowski, M. M., & Colangelo, N. (1984). Developmental potential of the gifted. *Gifted Child Quarterly, 28*(2), 80-88. https://doi.org/10.1177/001698628402800207

Qian, Y., Shuai, L., Chan, R. C., Qian, Q. J., & Wang, Y. (2013). The developmental trajectories of executive function of children and adolescents with attention deficit hyperactivity disorder. *Research in Developmental Disabilities, 34*(5), 1434-1445. https://doi.org/10.1016/j.ridd.2013.01.033

R Core Team. (2022). *R: A language and environment for statistical computing*. Vienna, Austria: R Foundation for Statistical Computing. Retrieved from http://www.R-project.org/

Rinn, A. N., & Reynolds, M. J. (2012). Overexcitabilities and ADHD in the gifted: An examination. *Roeper Review, 34*(1), 38-45. https://doi.org/10.1080/02783193.2012.627551

Roid, G. H., Miller, L. J., Pomplun, M., & Koch, C. (2013). *Leiter International Performance Scale-third edition*. Los Angeles, CA: Western Psychological Services.

Rommelse, N., Langerak, I., van der Meer, J., de Bruijn, Y., Staal, W., Oerlemans, A., & Buitelaar, J. (2015). Intelligence may moderate the cognitive profile of patients with ASD. *PLoS One, 10*(10), e0138698. https://doi.org/10.1371/journal.pone.0138698.

Ronald, A., de Bode, N., & Polderman, T. J. C. (2021). Systematic Review: How the Attention-Deficit/Hyperactivity Disorder Polygenic Risk Score Adds to Our Understanding of ADHD and Associated Traits. *Journal of the American Academy of Child & Adolescent Psychiatry, 60*, 1234-1277.